\documentclass[11pt]{article}
\usepackage{graphicx}

\def\card#1{\mathopen|#1\mathclose|}

\begin{document}

\title {Computational determination of (3,11) and (4,7) cages}

\author{Geoffrey Exoo \\
\small Department of Mathematics and Computer Science \\[-0.5ex]
\small Indiana State University \\[-0.5ex]
\small Terre Haute, IN 47809 \\
\tt gexoo@indstate.edu
\and
Brendan D. McKay\thanks{Corresponding author.  Supported by the
  Australian Research Council.} \\
\small School of Computer Science\\[-0.5ex]
\small Australian National University \\[-0.5ex]
\small Canberra, ACT 0200, Australia\\
\tt bdm@cs.anu.edu.au
\and
Wendy Myrvold and Jacqueline Nadon\\
\small Department of Computer Science\\[-0.5ex]
\small University of Victoria \\[-0.5ex]
\small Victoria, B.C., Canada V8W 3P6\\
\tt wendym@csc.uvic.ca
}

\maketitle

\begin{abstract}
A \textit{$(k,g)$-graph} is a $k$-regular graph of girth~$g$, and
a \textit{$(k,g)$-cage} is a $(k,g)$-graph of minimum order. We
show that a (3,11)-graph of order 112 found by Balaban in 1973 is
minimal and unique. We also show that the order of a (4,7)-cage
is~67 and find one example. Finally, we improve the lower bounds
on the orders of (3,13)-cages and (3,14)-cages to 202 and 260,
respectively. The methods used were a combination of heuristic
hill-climbing and an innovative backtrack search.
\end{abstract}

\begin{center}
AMS Subject Classifications: {05C25, 05C35}\\
Keywords: cage, regular graph, girth
\end{center}


\section{Introduction}

A {\em k-regular} graph is one in which every vertex has
degree $k$.
The {\em girth} of a graph is the length of a shortest
cycle.  
A {\em $(k,g)$-graph} is a regular graph of degree $k$ and
girth $g$, and
a {\em $(k,g)$-cage} is a $(k,g)$-graph of minimum possible order.

The problem of determining this order and identifying the
corresponding cages has been extensively studied. See~\cite{exoo}
for a current survey. The cases where the order of the cage was
known precisely before this work can be summarized as follows.

\begin{enumerate}
\itemsep=0pt
\item Regular graphs of degree 3 for girths up to 12;
\item Girth 5 graphs for degrees up to 7;
\item Girth 6, 8 and 12 graphs for degree one more than a prime power;
\item The case of degree 7 and girth 6.
\end{enumerate}

In this note, we add the case of degree 4 and girth 7 to the list,
and also show that the (3,11)-cage is unique.

The order, but not the uniqueness, of the $(3,11)$-cage was
previously announced in~\cite{soda}.

\section{Backtrack search}

The uniqueness of the (3,11)-cage, and the lower bound of 67 for
the (4,7)-cage, was proved using the program mentioned
in~\cite{soda} but not described in detail there.  We will
provide that description here using the (4,7)-cage
as an example.

\begin{figure}[ht]
\bigskip
\centering
\includegraphics[width=10cm]{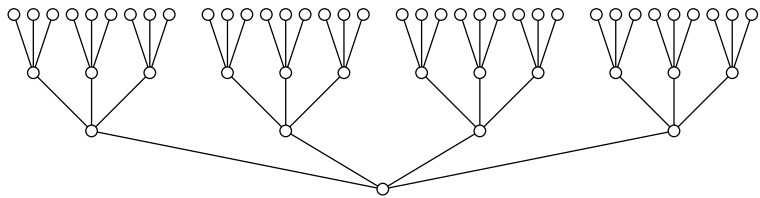}
\caption{The initial tree $T$}\label{tree47}
\end{figure}

Consider the construction of 4-regular graphs of girth at
least~7 and order $n\ge 53$. The vertices at distance at most~3
from some fixed vertex form a tree $T$ with 53 vertices as
in~Figure~\ref{tree47}. (In the case of even girth we would root the
tree at an edge rather than a vertex.)

Let $W$ be the set of $n-17$ vertices consisting of the 36 leaves
of $T$ and the $n-53$ vertices not in~$T$. The task is thus to add
additional edges within $W$ so that the resulting graph is quartic
and girth at least~$7$. We can do this in standard depth-first
manner, starting with the tree and adding one edge at a time. For
addition of a new edge $vw$ to be valid, $v$ and $w$ must have
degree less than~4 and be at distance at least~6; these properties
are easily monitored.

At each stage in the search, we choose one vertex of degree
less than~4 and try all the possibilities for joining it to
other vertices. We choose the vertex whose options are the most
restricted, as experiments showed this heuristic to be a good one.

\begin{figure}[ht]
\centering
\includegraphics[width=10cm]{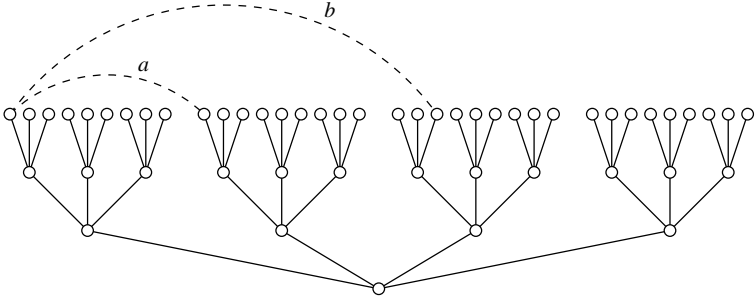}
\caption{Two equivalent choices}\label{tree47ab}
\end{figure}

Without additional improvements, this search is far too expensive
due to multiple equivalent subcases. The most obvious source of
equivalence is the automorphism group of the tree. Thus, the two
choices $a,b$ in Figure~\ref{tree47ab} are clearly equivalent and
there is no need to try $b$ as the first choice once $a$ has been
tried. The simple structure of the tree allows examples of this
type of equivalence to be monitored efficiently without explicit
computation of automorphism groups.

\begin{figure}[ht]
\centering
\includegraphics[width=10cm]{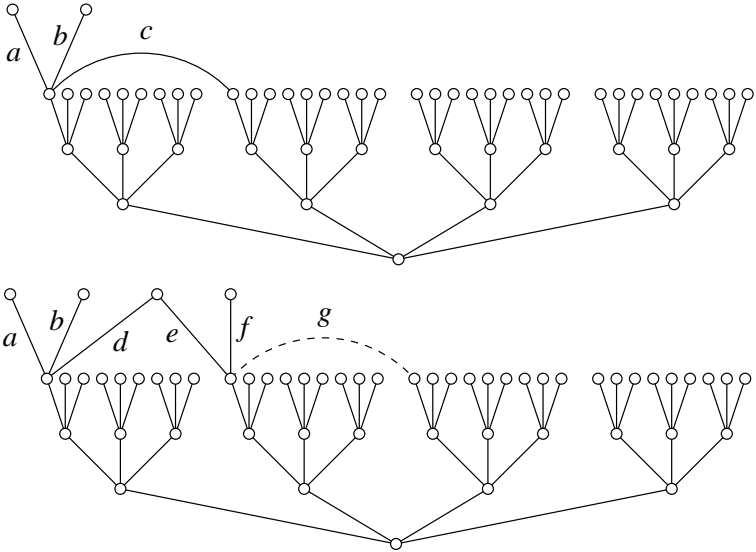}
\caption{Complex equivalence}\label{prune47c}
\end{figure}

However, when a moderate number of edges have been added,
discovery of automorphisms is considerably more difficult.
More importantly, there are equivalent subcases in the
search that do not derive directly from automorphisms of
any of the graphs that are constructed.  Consider
Figure~\ref{prune47c}.

Suppose we have completed the part of the search that begins with
the three edges $a,b,c$ in the upper diagram of the figure. This
means that, up to isomorphism, we have already found all the
quartic graphs that contain the upper diagram as a subgraph. Now
suppose that at a later point in the search we have added the edges
$a,b,d,e,f$ shown in the lower diagram of Figure~\ref{prune47c},
and are about to consider further edges. If we add edge $g$, then
we will not find any new quartic graphs, since $T+\{e,f,g\}$
is isomorphic to $T+\{a,b,c\}$. Therefore, we can mark $g$ as
ineligible. This applies whether we consider $g$ as the very next
edge to add or we consider adding it to some non-trivial extension
of $a,b,d,e,f$. For efficiency reasons we separate these two cases.

We now define the pruning process formally. Each node of the search
tree corresponds to a graph $T+E$, where $E$ is a set of edges
within~$W$. The \textit{subsearch rooted at $T+E$} is the subtree
of the search tree rooted at $T+E$, and $T+E$ is \textit{completed}
when we have finished the scan of that subtree.

\begin{description}
 \item[Pruning Rule. ] Suppose that $T+E$ is a search node and
    $vw$ is an available edge (that is, $v$ and $w$ have degree less
    than~4 and are at distance at least~6 from each other).
    Suppose further that there
    is a subset $E'\subseteq E+vw$ such that $T+E'$ is
    isomorphic to a node which is already completed.
    Then the subsearch rooted at $T+E$ can avoid adding
    the edge~$vw$.
\end{description}

The pruning rule is applied using the \texttt{nauty} graph
isomorphism software~\cite{nauty}. Since subgraph isomorphism
testing is required, the rule is very expensive compared to the time
otherwise required to process one node of the search. In practice
we limit it to only certain $E$ and certain $v,w$. Specifically,
we define two levels $\ell_1\ge\ell_2\ge0$. For search nodes $T+E$
such that $\card{E}\le\ell_2$, we apply the rule for all~$v,w$. For
$\ell_2<\card{E}\le\ell_1$, we apply it only to those $v,w$ such
that $vw$ is a candidate for the very next edge to add.

For larger $\ell_1,\ell_2$, the number of nodes in the search tree
is reduced but the time expended in pruning the tree is increased.
So there is some optimal compromise.

A further technique is needed in very difficult cases (such as
$n=66$), for which it is desirable to divide the computation across
multiple processes (perhaps on different computers). The division
is accomplished in the usual fashion: for some level $\ell_0$, treat
the subsearchs rooted at $T+E$ for $\card{E}=\ell_0$ as independent.
Subsearch $i$ is assigned to process $j$ if $i\equiv j~(\textrm{mod}\;N)$,
where $N$ is the number of processes.

Each process involved in the search computes the whole search
tree to level $\ell_0$ and only its own subsearches at higher
levels. Since the application of the pruning rule is hard to achieve
across multiple processes, we use $\ell_0>\ell_1$. Higher $\ell_0$
also tends to share the load more evenly between processes. These
considerations, however, imply that a large part of the tree,
including all the rule applications, are repeated in every process,
which would seem to prevent use of large $\ell_1,\ell_2$ due to the
expensive nature of the pruning rule. The solution is to conduct the
computation in two phases. In the first phase, executed on a single
processor, the search tree is computed up to level $\ell_1$ while
applying the pruning rule. During this phase, the results of all
the rule applications are recorded in an audit file. Then, in the
second phase where the full search tree is computed, the pruning
rule is applied at almost no cost by following the audit file. If
$\ell_0$ is not too high, the cost of the division of the second
phase into parts is negligible.

For $n=66$, we used $\ell_0=22$, $\ell_1=16$ and $\ell_2=6$.
The total number of nodes in the search tree was
318,904,129,273,923.  Using a large mix of computers
ranging from 650\,MHz to 3\,GHz, the total time was
96~years (105,000 nodes per second).
This cost was about 5 times the cost for $n=65$.

Improvements to the program after the computation finished resulted
in a 4-fold speedup. Nevertheless, the computation of all the cages
of order~67 will not be feasible in the near future.

Exhaustive computation of the $(3,11)$-cages, showing that the
graph found by Balaban~\cite{balaban} is the only one, required
17 years on computers averaging about 300 MHz.

With lesser expenditure of cpu time, we also showed that a
$(3,13)$-cage has at least 202 vertices (improved from 196),
and that a $(3,14)$-cage has at least 260 vertices
(improved from 256). 
The first result, and the bound 258 for the second,
were previously announced in~\cite{soda}.

To emphasise the efficiency of this method, we compared it against
the previously best code for cubic graphs~\cite{maus}. For girth 9
and 58 vertices, our approach was 52 times faster, while for girth
11 and 104 vertices our approach was 856 times faster. All of the
results we highlight in this paper were previously infeasible.

\section{A Graph of Degree 4 and Girth 7}

As noted in the previous section, we showed exhaustively that
there are no $(4,7)$-graphs with 66 or fewer vertices.  To
complete the proof that $(4,7)$-cages have 67 vertices, we
present an example.

Our example on 67 vertices is given below as an adjacency list. The
graph has an automorphism group of order four, isomorphic to $Z_2
\times Z_2$. In the action of this group on the vertices, there are
11 orbits of length 4, 11 orbits of length 2, and one fixed point.

The graph was constructed using a hill-climbing algorithm that
begins with an empty graph on 67 vertices and adds edges, one at
a time, while not violating the degree and girth conditions. An
outline of the algorithm follows.

\begin{tt}
\begin{tabbing}
while \= there are vertices with degree < 4: \\
\>  while \= there are edges that can be added: \\
\>  \> for \= each edge that can be added \\
\>  \>   \> compute the degree sum of its vertices \\
\>  \> pick the edge with the largest sum \\
\>  \> choose randomly in case of ties \\
\>  \> add the winning edge \\
\>  if all vertices have degree 4: \\
\>  \> save graph \\
\>  \> exit \\
\>  delete 1, 2 or 3 random edges
\end{tabbing}
\end{tt}
 
\noindent A few of the steps require some explanation.

\begin{enumerate}
\item Edges are added
between vertices with maximum possible degree sum.
Doing this in the early stages of the algorithm appears to be
essential to finding a cage.  We ran our program many times with
this condition removed, without success.
\item When edges are added, a record is made of when they
were added.  Time is measured in trips through the outer 
loop.
\item When edges are chosen for deletion, it is done in one
of two ways.  Either the probability that
an edge is chosen is proportional to its age, or it is
inversely proportional to its age.  These two modes are
alternated, each being used for a few thousand time periods.
Very recently deleted edges are not chosen for reinsertion.
\end{enumerate}

A program implementing this method has been run several hundred
times. It succeeds in finding a graph approximately twenty percent
of the time. Each time it has succeeded, it has found the same
graph.

\begin{figure}
\centering
\def\0{\hphantom{0}}
\setlength{\tabcolsep}{2em}
\def\arraystretch{0.8}
\begin{tabular}{lll}
\00: \01 \02 \03 \04 & 23: \07 61 43 29 & 46: 14 31 41 24 \\
\01: \00 \05 \06 \07 & 24: \07 26 36 46 & 47: 15 42 58 26 \\
\02: \00 \08 \09 10  & 25: \07 63 40 34 & 48: 15 34 35 19 \\
\03: \00 11 12 13    & 26: \08 53 24 47 & 49: 15 39 65 29 \\
\04: \00 14 15 16    & 27: \08 43 45 22 & 50: 16 33 66 43 \\
\05: \01 17 18 19    & 28: \08 52 18 56 & 51: 16 36 54 30 \\
\06: \01 20 21 22    & 29: \09 23 49 59 & 52: 16 28 21 40 \\
\07: \01 23 24 25    & 30: \09 19 62 51 & 53: 59 66 38 26 \\
\08: \02 26 27 28    & 31: \09 20 46 57 & 54: 58 64 51 22 \\
\09: \02 29 30 31    & 32: 10 44 65 21  & 55: 39 62 33 56 \\
10: \02 32 33 34     & 33: 10 50 55 17  & 56: 28 55 35 61 \\
11: \03 35 36 37     & 34: 10 48 64 25  & 57: 58 31 17 40 \\
12: \03 38 39 40     & 35: 11 56 48 20  & 58: 54 57 47 61 \\
13: \03 41 42 43     & 36: 11 51 65 24  & 59: 64 37 29 53 \\
14: \04 44 45 46     & 37: 11 59 17 45  & 60: 18 66 63 65 \\
15: \04 47 48 49     & 38: 12 53 19 44  & 61: 58 44 23 56 \\
16: \04 50 51 52     & 39: 12 55 49 22  & 62: 55 30 63 42 \\
17: \05 57 37 33     & 40: 12 52 25 57  & 63: 25 62 60 45 \\
18: \05 60 28 41     & 41: 13 46 18 64  & 64: 59 54 34 41 \\
19: \05 38 30 48     & 42: 13 62 47 21  & 65: 36 32 49 60 \\
20: \06 35 31 66     & 43: 13 23 27 50  & 66: 50 60 53 20 \\
21: \06 32 52 42     & 44: 14 32 38 61  & \\
22: \06 54 39 27     & 45: 14 27 37 63  &
\end{tabular}
\caption{Adjacency list for a (4,7)-graph of order 67.}
\end{figure}


\end{document}